

\baselineskip=14pt
\parskip=10pt

\font\eightrm=cmr8 

\magnification=\magstephalf

\def\1{{\overline{1}}}
\def\2{{\overline{2}}}
\parindent=0pt
\overfullrule=0in

\def\frac#1#2{{#1 \over #2}}
\centerline
{\bf Lots and Lots of Perrin-Type Primality Tests and Their Pseudo-Primes }
\bigskip
\centerline
{\it Robert DOUGHERTY-BLISS and Doron ZEILBERGER}

{\bf Abstract:}  We use {\it Experimental Mathematics} and {\it Symbolic Computation} (with Maple), to search for lots and lots
of Perrin- and Lucas- style primality tests, and try to sort the wheat from the chaff. More impressively, we find quite
a few such primality tests for which we can explicitly construct infinite families of pseudo-primes, rather, like in the
cases of Perrin pseudo-primes and the famous Carmichael primes, proving the mere existence of infinitely many of them.

{\bf Preface: How it all Started thanks to Vince Vatter}

It all started when we came across Vince Vatter's delightful article [V], where he gave a cute combinatorial proof, inspired by COVID, and social distancing,
of the following fact that goes back to Raoul Perrin [P] (See also [Sl1], [Sl2], [St], [W]).

{\bf Perrin's Observation}:  Let the integer sequence $A(n)$ defined by:
$$
A(1)=0 \quad,\quad
A(2)=2 \quad,\quad
A(3)=3 \quad,\quad
A(n)=A(n-2)+A(n-3) \,\,( for \quad n>3) \quad,
$$
then for every prime $p$, we have:
$$
p|A(p) \quad .
$$

Perrin, back in 1889, was wondering whether the condition is {\bf sufficient}, i.e. whether there are any {\it pseudo-primes}, i.e. {\it composite} $n$ such that
$A(n)/n$ is an integer. He could not find any, and as late as 1981, none was found $\leq 140000$ (see [AS]). In 1982, Adams and Shanks [AS] {\it rather quickly}
found the smallest  Perrin pseudo prime, $271441$, followed by the next-smallest, $904631$, and then they found quite a few other ones. 
Jon Grantham [Gr] proved that there are {\it infinitely many} Perrin pseudo-primes, and finding as many as possible of them, became a computational challenge, see Holger's paper [H].

Another, older, primaility test is that based on the Lucas numbers ([Sl3], [Sl4]).

{\bf Vince Vatter's Combinatorial Proof}

Vatter first found a {\bf combinatorial interpretation} of the Perrin numbers, as the number of circular words of length $n$ in the alphabet $\{0,1\}$, that {\bf avoid} the 
{\bf consecutive subwords} (aka {\it factors} in formal language lingo), $\{000, 11\}$.

More formally: words $w=w_1, \dots, w_n$ in the alphabet $\{0,1\}$, such that for $1 \leq i \leq n-2$, $w_iw_{i+1}w_{i+2} \neq 000$, and also
$w_{n-1}w_n w_1 \neq 000$ and $w_{n}w_1 w_2 \neq 000$ as well as
for $1 \leq i \leq n-1$, $w_iw_{i+1} \neq 11$, and $w_nw_1 \neq 11$.

Then he argued that if $p$ is a prime, all the $p$ circular shifts are {\bf different}, since otherwise there would be a non-trivial period, that can't happen since $p$ is prime. 
Since the constant words $0^p$ and $1^p$ obviously can't avoid both $00$ and $111$, Perrin's theorem follows.

This proof is reminiscent of Solomon Golomb's [G] snappy combinatorial proof of Fermat's little theorem [G] that argued that there are
$a^p-a$ non-monochromatic straight necklaces with $p$ beads of $a$ colors, and for each such necklace, the $p$ rotations are all different (see also [Z1], p. 560).

When we saw Vatter's proof we got all excited. Vatter's argument transforms {\it verbatim} to counting circular words in {\it any} (finite) alphabet, and any (finite) set of
forbidden (consecutive) patterns! More than twenty years ago one of us (DZ) wrote a paper, in collaboration with his then PhD student, Anne Edlin [EZ], that {\it automatically}
finds the (rational) generating function in any such scenario, hence this is a cheap way to manufacture lots and lots of Perrin-style primality tests.
We already had a  Maple package {\tt https://sites.math.rutgers.edu/\~{}zeilberg/tokhniot/CGJ} to handle it, so all that remained was to {\it experiment} with many alphabets
and many sets of forbidden patterns, and search for those that  have only few small {\it pseudo-primes}.

This inspired us to write our first Maple package, {\tt PerrinVV.txt}, available from

{\tt https://sites.math.rutgers.edu/\~{}zeilberg/tokhniot/PerrinVV.txt}  \quad .

See the front of this article

{\tt https://sites.math.rutgers.edu/\~{}zeilberg/mamarim/mamarimhtml/perrin.html} \quad,

for many such primality tests, inspired by sets of forbidden patterns, along with all the pseudoprimes less than a million.

{\bf An even better way to manufacture Perrin-style Primality tests}

After the initial excitement we got an {\it epiphany}, and as it turned out, it was already made, in 1990, by Stanley Gurak [Gu]. Take {\it any}  polynomial
$Q(x)$ with {\bf integer coefficients}, and constant term $1$, and write it as
$$
Q(x)=1 \,- \, e_1\,x \, +\, e_2x^2 - \dots +(-1)^k\, e_k x^k \quad.
$$
Factorize it over the complex numbers 
$$
Q(x) \,= \, (1- \alpha_1 x )  (1- \alpha_2 x )  \cdots (1- \alpha_k x )  \quad .
$$
Note that $e_1,e_2, \dots$ are the {\bf elementary symmetric functions} in $\alpha_1, \dots, \alpha_k$.

Defining
$$
a(n):= \alpha_1^n \, + \, \alpha_2^n \, + \, \dots \, + \alpha_k^n \quad,
$$
it follows thanks to Newton's identities ([M]) that $\{a(n)\}$ is an {\bf integer sequence}. The generating function
$$
\sum_{n=0}^{\infty} a(n)\,x^n \, = \, \frac{1}{1-\alpha_1x} \, + \, \frac{1}{1-\alpha_2x} \, + \, \dots \, + \,   \frac{1}{1-\alpha_kx} \quad ,
$$
has denominator $Q(x)$ and {\it some} numerator, let's call it $P(x)$, with {\bf integer coefficients}, that Maple can easily find {\it all by itself}.

So we can {\it define} an integer sequence $\{a(n)\}$ in terms of the rational function $P(x)/Q(x)$, where $Q(x)$ is {\it any} polynomial with constant term $1$, and
$P(x)$ {\it comes out} as above:

$$
\sum_{n=0}^{\infty} a(n) x^n \, = \, \frac{P(x)}{Q(x)} \quad .
$$

We claim that each such integer sequence engenders a {\it Perrin-style}
primality test, namely
$$
a(p) \equiv e_1 (mod \,\, p).
$$
To see this, note that
$$
(\alpha_1 + \cdots + \alpha_k)^p = a(p) \,+\, p\,A(p),
$$
where
$$
A(p) = \sum_{
{{i_1 + i_2 + \cdots + i_k = p} \atop {i_1, i_2, \dots i_k < p}}
} \frac{(p-1)!}{i_1! \cdots i_k!} \alpha_1^{i_1} \cdots \alpha_k^{i_k} 
$$
is a symmetric polynomial in the $\alpha_i$. The fundamental theorem of
symmetric functions [M] implies that $A(p)$ is an integer. 
Fermat's little theorem then gives
$$
a(p) \equiv (\alpha_1 + \cdots + \alpha_k)^p = e_1^p \equiv e_1 (mod \,\, p).
$$

So this is an even easier way to manufacture lots and lots of Perrin-style primality tests, and we can let the computer search for those
that have as few small pseudo-primes as possible.

This is implemented in the Maple package {\tt Perrin.txt},  available from

{\tt https://sites.math.rutgers.edu/\~{}zeilberg/tokhniot/Perrin.txt}  \quad .

Again, see the front of this article

{\tt https://sites.math.rutgers.edu/\~{}zeilberg/mamarim/mamarimhtml/perrin.html} \quad,

for many such primality tests, inspired by this more general method (first suggested by Stanley Gurak [Gu]).

Searching for such primality tests with as few pseudo-primes less than a million, we stumbled on the following:

{\bf The DB-Z Primality Test}

Let 
$$
\sum_{n=0}^{\infty} a(n)\, x^n \, := \, \frac{-3 x^{4}-5 x^{2}-6 x +7}{-4 x^{7}-x^{4}-x^{2}-x +1} \quad,
$$
or equivalently, the integer sequence defined by
$$
a(1)=1 \,,\,  a(2)=3 \,,  \, a(3)=4 \,, \, a(4)=11 \, , \, a(5)=16 \,, \,a(6)=30 \,, \,a(7)=78,
$$
$$
a(n)\,=\,a(n-1)\,+\,a(n-2) \,+ \, a(n-4) \,+ \, 4\,a(n-7) \,\,( for \quad n>7) \quad,
$$
then if $p$ is prime, we have
$$
a(p) \equiv 1\, (mod \,\, p) \quad.
$$

Manuel Kauers kindly informed us that the seven smallest DB-Z pseudo-primes are

$\bullet$ $1531398 = 2 \cdot 3 \cdot 11 \cdot 23203$ \quad,
    
$\bullet$  $114009582 = 2 \cdot 3 \cdot 17 \cdot 1117741$ \quad,
    
$\bullet$  $940084647 = 3 \cdot 47 \cdot 643 \cdot 10369$ \quad,
    
$\bullet$  $4206644978 = 2 \cdot 97 \cdot 859 \cdot 25243$ \quad,
    
$\bullet$  $7962908038 = 2 \cdot 191 \cdot 709 \cdot 29401$ \quad,

$\bullet$ $20293639091 = 11 \cdot  3547 \cdot 520123$ \quad,

$\bullet$ $41947594698= 2\cdot 3 \cdot 19 \cdot 523 \cdot 703559 $ \quad.

{\eightrm [This was a computational challenge posed by us to Manuel Kauers, and we offered to donate $100$ dollars to the OEIS in his honor. The donation was made].}

{\bf An even better  Primality Test}

After the first version of this paper was written, with the help of Manuel Kauers, we discovered an even better primality test.

{\bf The DB-Kauers primality test}

Let 
$$
\sum_{n=0}^{\infty} a(n)\, x^n \, := \, \frac{-9 x^{5}-16 x^{4}-10 x +6}{-3 x^{6}-9 x^{5}-8 x^{4}-2 x +1}
\quad,
$$
or equivalently, the integer sequence defined by
$$
a(1)=2 \,,\,  a(2)=4 \,,  \, a(3)=8 \,, \, a(4)=48 \, , \, a(5)=157 \,, \,a(6)=382 \,,
$$
$$
a(n)\,=\,2a(n-1)\,+ \, 8\, a(n-4) \,+ \, 9\,a(n-5) \,+\, a(n-6)\, \,( for \quad n>6) \quad,
$$
then if $p$ is prime, we have
$$
a(p) \equiv 2\, (mod \,\, p) \quad.
$$

The smallest pseudoprime happens to be  $2,260,550,373 = 3 \cdot 103 \cdot 107 \cdot 68371$.

{\bf Perrin-Style Primality Tests with Explicit Infinite Families of Pseudo-Primes}

We are particularly proud of the next primality test, featuring  the  {\it Companion Pell numbers} {\tt https://oeis.org/A002203}. These numbers 
have been studied extensively, but as far as we know using them as a {\it primality test} is new. It is not a  very good one, but
the novelty is that it has an {\it explicit} {\bf doubly-infinite} set of pseudo-primes.

{\bf The Companion Pell Numbers Primality Test}
Let 
$$
\sum_{n=0}^{\infty} a(n)\, x^n \, := \, \frac{2-2 x}{-x^{2}-2 x +1} \quad,
$$
or equivalently,
$$
a(1)=2 \,,\, a(2)=6 \quad \,, \,\quad a(n)\,=\,a(n-1)\,+ \, 2a(n-2) \quad (for \quad n>2) \quad,
$$
then if $p$ is a prime, we have
$$
a(p) \equiv 2 \,(mod \,\, p) \quad .
$$

{\bf Theorem 1}; The following doubly-infinite family
$$
\{\,2^i \cdot 3^j \, | \, i \geq 3 \, ,  \, j \geq 0 \,\} \quad,
$$
are Companion-Pell pseudoprimes, in other words,
$$
\frac{a(2^i \cdot 3^j)-2}{2^i 3^j} \quad,
$$
are always integers, if $i \geq 3$ and $j \geq 0$.

{\bf Proof}: let
$$
\alpha_1:=1+\sqrt{2} \quad, \quad \alpha_2:=1-\sqrt{2} \quad, 
$$
then, since
$$
\frac{2-2 x}{-x^{2}-2 x +1} \, = \, \frac{1}{1-\alpha_1\,x} \, + \, \frac{1}{1-\alpha_2\,x} \quad,
$$
we  have the Binet-style formula
$$
a(n)=\alpha_1^n+ \alpha_2^n \quad.
$$
Since $\alpha_1 \cdot \alpha_2=-1$, we have the following two recurrences (check!)
$$
a(2n)=a(n)^2 +2 (-1)^{n+1} \quad,
$$
$$
a(3n)=a(n)^3 +3 (-1)^{n+1} a(n) \quad.
$$

Let 
$$
b(n)=a(n)-2 \quad,
$$
hence we have

{\bf Fact 1}: If $n$ is even then
$$
b(2n)=b(n)(b(n)+4) \quad .
$$

It follows immediately that:

{\it if $n$ is even and $b(n)/n$ is an integer then $b(2n)/(2n)$ is also an integer.}

It remains to prove that $b(8\,3^j)/(8 \, 3^j)$ is an integer, for all $j\geq 0$.

Thanks to the second recurrence we have

$$
b(3n)+2=(b(n)+2)^3 +3 (-1)^{n+1} (b(n)+2) \quad.
$$

Hence

{\bf Fact 2}: If $n$ is even then
$$
b(3n)=b(n)(b(n)^2 \,+ \, 6\,b(n) \,+ \,12) \quad .
$$

It follows immediately that:

{\it if $n$ is divisible by $6$ and $b(n)/n$ is an integer then $b(3n)/(3n)$ is also an integer.}

Since $b(24)/24$ is an integer, the theorem follows by induction.

We now state without proofs (except for Theorem 4, where we give a sketch) a few other primality tests that have explicit infinite families of  pseudoprimes.

{\bf Theorem 2}: Let
$$
\sum_{n=0}^{\infty} a(n)\, x^n \, := \, \frac{2-x}{-2 x^{2}-x +1} \quad,
$$
or equivalently,
$$
a(1)=1 \, , \, a(2)=5 \quad , \quad a(n)\,=\,a(n-1) \,+ \, 2a(n-2) \quad (for \quad n>2) \quad,
$$
then if $p$ is a prime, we have
$$
a(p) \equiv 1 \,(mod \,\, p) \quad .
$$
Furthermore, $\{2^i \,|  i \geq 2\}$ are all pseudo-primes, in other words
$$
a(2^i) \equiv 1 \,(mod \,\, 2^i) \quad , \quad i \geq 2.
$$

{\bf Theorem 3}: Let
$$
\sum_{n=0}^{\infty} a(n)\, x^n \, := \, \frac{2-2 x}{-2 x^{2}-2 x +1} \quad,
$$
or equivalently,
$$
a(1)=2 \,, \,  a(2)=8 \quad , \quad a(n)\,=\, 2a(n-1) \,+ \, 2a(n-2) \quad (for \quad n>2) \quad,
$$
then if $p$ is a prime, we have
$$
a(p) \equiv 2 \,(mod \,\, p) \quad .
$$
Furthermore, the following infinite families are all pseudo-primes:

$\{3^i \,|  i \geq 2\}$ \quad, \quad  $\{2\cdot 3^i \,|  i \geq 1\}$ \quad , \quad  $\{11\cdot 81^i \,|  i \geq 1\}$, $\{23\cdot 3^{5i} \,|  i \geq 1\}$ \quad, \quad
$\{29\cdot 3^{4+12i} \,|  i \geq 0\}$ \quad, \quad $\{31\cdot 3^{16i} \,|  i \geq 1\}$ .

{\bf Theorem 4}: Let
$$
\sum_{n=0}^{\infty} a(n)\, x^n \, := \, \frac{2 x^{2}+3}{2 x^{3}+2 x^{2}+1} \quad,
$$
or equivalently,
$$
a(1)=0 \, , \,  a(2)=-4 \, , \, a(3)=-6\quad , \quad a(n)\,=\, -2a(n-2) \,- \, 2a(n-3) \quad (for \quad n>2) \quad,
$$
then if $p$ is a prime, we have
$$
a(p) \equiv 0 \,(mod \,\, p) \quad .
$$
Furthermore, the following infinite families are all pseudo-primes:

$\{2^i \,|  i \geq 2\}$ \quad ,  \quad $\{3\cdot2^{4i} \,|  i \geq 2\}$ \quad,  \quad  $\{11\cdot2^{18i} \,|  i \geq 2\}$ \quad ,  \quad $\{13\cdot2^{17+20i} \,|  i \geq 2\}$ .

{\bf Sketch of Proof}: We use the {\it C-finite ansatz} [Z2]. Let
$$
b(n)=a(2n)-a(n)^2 \quad,
$$
then it follows from the C-finite anzatz that $b(n)$ satisfies {\it some} recurrence,  that turns out to be
$$
b(1)=-4 \, , \, b(2)=-8 \, , \, b(3)= -40 \quad, \quad b(n)\,=\, 2b(n-1)+ 4 b(n-3) \,\, \,(for \quad n>3)
$$
We now define
$$
c(n):=\frac{b(n)}{2^{\lfloor n/2 \rfloor}} \quad,
$$
and once again it follows that $c(n)$ satisfies the recurrence, 
$$
c(1)=-4 \, , \, c(2)=-4 \, , \, c(3)= -20 , \, c(4)= -24 , \, c(5)= -56 , \, c(6)= -76 
$$
$$
c(n)\,=\, 2c(n-2)+ 4 c(n-4)+ 2 c(n-6) \,\,\,(for \quad n>6) \quad.
$$
Note that $c(n)$ are manifestly {\bf integers}. Going back to $a(n)$ we have the recurrence
$$
a(2n)=a(n)^2+ 2^{\lfloor n/2 \rfloor} c(n) \quad,
$$
and it follows by induction that $a(2^i)/2^i$ are all integers. A similar argument goes for the other infinite families claimed.

{\bf Theorem 5}: Let
$$
\sum_{n=0}^{\infty} a(n)\, x^n \, := \, \frac{-2 x^{2}-2 x +3}{-x^{3}-2 x^{2}-x +1} \quad,
$$
or equivalently,
$$
a(1)=1 \, , \,  a(2)=5 \, , \, a(3)=10 \quad , \quad a(n)\,=\, a(n-1) \,+ \, 2a(n-2) +a(n-3)\quad (for \quad n>2) \quad,
$$
then if $p$ is a prime, we have
$$
a(p) \equiv 1 \,(mod \,\, p) \quad .
$$
Furthermore, the following infinite families are all pseudo-primes:

$\{3^i \,|  i \geq 2\}$  \quad ,  \quad  $\{5 \cdot 3^{6+10\,i} \,|  i \geq 0\}$ \quad ,  \quad  $\{5 \cdot 3^{8+10\,i} \,|  i \geq 0\}$  \quad ,  \quad  $\{7 \cdot 3^{4+6\,i} \,|  i \geq 0\}$ ,

We found $9$ other such primality tests, with infinite explicit families of presodoprimes, that can  be viewed by typing

{\tt PDB(x);} \quad ,

in the Maple package {\tt Perrin.txt}.

For fast computations and explorations using $C$ programs, readers are welcome to explore
RDB's github site:

{\tt https://github.com/rwbogl/pt } \quad .

{\bf Acknowledgment}:  Many thanks to Manuel Kauers for his computational prowess, and to the referee for a helpful remark.

{\bf References}

[AS] William Adams and Daniel Shanks, {\it Strong primality tests that are not sufficient}, Mathematics of Computation {\bf 39} (1982), 255-300.

[EZ] Anne E. Edlin and Doron Zeilberger, {\it  The Goulden-Jackson Cluster method For cyclic  Words}, Advances in Applied Mathematics {\bf 25}(2000), 228-232.  \hfill\break
{\tt https://sites.math.rutgers.edu/\~{}zeilberg/mamarim/mamarimhtml/cgj.html} \quad .

[Go] Solomon W. Golomb, {\it Combinatorial Proof of Fermat's ``Little'' Theorem}, The American Mathematical Monthly, {\bf 63}(10),  (Dec., 1956), 718. \hfill\break
{\tt https://sites.math.rutgers.edu/\~{}zeilberg/akherim/golomb56.pdf} \quad .

[Gr] Jon Grantham, {\it There are infinitely many Perrin pseudoprimes}, Journal of Number Theory {\bf 130} (2010), 1117-1128.

[Gu] Stanley Gurak, {\it Pseudoprimes for Higher-Order Linear Recurrence Sequences}, 
Mathematics of Computation, {\bf 55} (1990), 783-813.

[H] Stephan Holger, {\it Millions of Perrin pseudoprimes including a few giants}, arXiv:2002.03756 [math.NA], 2020. {\tt https://arxiv.org/abs/2002.03756} \quad.

[M] Ian G. Macdonald, {\it ``Symmetric Functions and Hall Polynomials''}, Second Edition, Clarendon Press, Oxford, 1995.

[P] Raoul Perrin, {\it Item 1484}, L'Interm\`ediare des math {\bf 6} (1899), 76-77.

[Sl1] Neil J. A. Sloane, {\it The On-line Sequence of Integer Sequences}, Sequence  A001608, \hfill\break
{\tt https://oeis.org/A001608} \quad. 

[Sl2] Neil J. A. Sloane, {\it The On-line Sequence of Integer Sequences}, Sequence  A013998, \hfill\break   
{\tt  https://oeis.org/A013998} 

[Sl3] Neil J. A. Sloane, {\it The On-line Sequence of Integer Sequences}, Sequence  A005854, \hfill\break   
{\tt  https://oeis.org/A005845} \quad. 

[Sl4] Neil J. A. Sloane, {\it The On-line Sequence of Integer Sequences}, Sequence  A000032, \hfill\break   
{\tt   https://oeis.org/A000032} \quad. 

[St] Ian Stewart, {\it Tales of a Neglected Number}, Mathematical Recreations, Scientific American {\bf 274}(6) (1996), pp. 102-103.

[V] Vince Vatter, {\it Social Distancing, Primes, and Perrin Numbers}, , Math Horizons, {\bf 29}(1) (2022). \hfill\break
{\tt https://sites.math.rutgers.edu/\~{}zeilberg/akherim/vatter23.pdf} \quad .

[W] Wikipedia, {\it Perrin Number}, {\tt https://en.wikipedia.org/wiki/Perrin$\_$number} \quad .

[Z1] Doron Zeilberger, {\it Enumerative and Algebraic Combinatorics}, in: Princeton Companion to Mathematics (edited by W.T. Gowers),
Princeton University Press, 2008, 550-561. \hfill\break
{\tt https://sites.math.rutgers.edu/\~{}zeilberg/mamarim/mamarimPDF/enu.pdf} \quad .

[Z2] Doron Zeilberger, {\it The C-finite ansatz}, Ramanujan Journal {\bf 31} (2013), 23-32. \hfill\break
{\tt https://sites.math.rutgers.edu/\~{}zeilberg/mamarim/mamarimhtml/cfinite.html} \quad .
\bigskip
\hrule
\bigskip
Robert Dougherty-Bliss , Department of Mathematics, Rutgers University (New Brunswick), Hill Center-Busch Campus, 110 Frelinghuysen
Rd., Piscataway, NJ 08854-8019, USA. \hfill\break
Email: {\tt  robert dot w dot bliss at gmail dot com}   \quad .
\bigskip
Doron Zeilberger, Department of Mathematics, Rutgers University (New Brunswick), Hill Center-Busch Campus, 110 Frelinghuysen
Rd., Piscataway, NJ 08854-8019, USA. \hfill\break
Email: {\tt DoronZeil at gmail  dot com}   \quad .

\bigskip
Published in INTEGERS {\bf 23} (2023), \#A95.
\bigskip
First written: July 15, 2023; This version (correcting typos): March 30, 2024.

\end